\documentclass[twoside,leqno, british,12pt]
{amsart}
\usepackage{amsmath, amssymb, amsbsy, amsthm}
\usepackage{eucal,url,amssymb,
enumerate,amscd,stmaryrd
}

\usepackage[T1]{fontenc}

\usepackage[british]{babel}
\usepackage[pagebackref,colorlinks=true]{hyperref}  

\begin{document}

\begin{abstract}
The object of investigations are almost contact B-metric manifolds which are derived as a product of a real line and a 2-dimensional manifold equipped with a complex structure and a Norden metric. There are used two different methods for generation of the B-metric on the product manifold. The constructed manifolds are characterised with respect to the Ganchev-Mihova-Gribachev classification and their basic curvature properties.
\end{abstract}


\title[Almost contact B-metric manifolds as extensions]
{Almost contact B-metric manifolds as extensions of a 2-dimensional space-form}

\author[H. Manev]{Hristo Manev$^{1,2}$}
\address[1]{Medical University of Plovdiv, Faculty of Pharmacy,
Department of Pharmaceutical Sciences,   15-A Vasil Aprilov
Blvd.,   Plovdiv 4002,   Bulgaria;}
\address[2]{Paisii Hilendarski University of Plovdiv,   Faculty of Mathematics and
Informatics,   Department of Algebra and Geometry,   236
Bulgaria Blvd.,   Plovdiv 4027,   Bulgaria}
\email{hmanev@uni-plovdiv.bg}


\frenchspacing

\newcommand{\ie}{i.e. }
\newcommand{\X}{\mathfrak{X}}
\newcommand{\W}{\mathcal{W}}
\newcommand{\LL}{\mathcal{L}}
\newcommand{\F}{\mathcal{F}}
\newcommand{\M}{(M,\allowbreak{}\f,\allowbreak{}\xi,\allowbreak{}\eta,\allowbreak{}g)}
\newcommand{\R}{\mathbb{R}}
\newcommand{\C}{\mathbb{C}}
\newcommand{\n}{\nabla}
\newcommand{\al}{\alpha}
\newcommand{\f}{\varphi}
\newcommand{\ta}{\theta}
\newcommand{\lm}{\lambda}
\newcommand{\om}{\omega}
\newcommand{\norm}[1]{\Vert#1\Vert ^2}
\newcommand{\Span}{\mathrm{span}}
\newcommand{\Id}{\mathrm{Id}}
\newcommand{\D}{\mathrm{d}}
\newcommand{\ddr}{\tfrac{\D}{\D r}}
\newcommand{\g}{\check{g}}
\newcommand{\nn}{\check{\n}}
\newcommand{\ddu}[1]{\frac{\partial}{\partial u^{#1}}}
\newcommand{\ddv}[1]{\frac{\partial}{\partial v^{#1}}}
\newcommand{\ddt}{\tfrac{\D}{\D t}}

\newcommand{\thmref}[1]{The\-o\-rem~\ref{#1}}
\newcommand{\propref}[1]{Pro\-po\-si\-ti\-on~\ref{#1}}
\newcommand{\secref}[1]{\S\ref{#1}}
\newcommand{\lemref}[1]{Lem\-ma~\ref{#1}}
\newcommand{\dfnref}[1]{De\-fi\-ni\-ti\-on~\ref{#1}}
\newcommand{\corref}[1]{Corollary~\ref{#1}}



\newtheorem{thm}{Theorem}[section]
\newtheorem{lem}[thm]{Lemma}
\newtheorem{prop}[thm]{Proposition}
\newtheorem{cor}[thm]{Corollary}
\newtheorem{conv}[thm]{Convention}
\newtheorem{rmk}[thm]{Remark}

\theoremstyle{definition}
\newtheorem{defn}{Definition}[section]

\hyphenation{Her-mi-ti-an ma-ni-fold ah-ler-ian}




\maketitle
%
%


\section{Introduction}
The differential geometry of almost contact metric manifolds
is well studied (e.g. \cite{Blair}). Ganchev, Mihova, Gribachev begin
investigations on the almost contact manifolds with B-metric in \cite{GaMiGr}.
These manifolds are the odd-dimensional
counterpart of almost complex manifolds with Norden metric (briefly, almost Norden manifold{s})
\cite{GaBo,GrMeDj}, where the almost complex structure acts as an
anti-isometry regarding the metric. Further, almost contact B-metric manifolds {of arbitrary odd dimension} are
investigated in many works, for example \cite{Man31,ManGri1,ManGri2,ManIv13,ManIv14,NakGri2}.

An object of our special interest is the case of the lowest dimension 3 of the considered manifolds (\cite{HM1,HM4,HM5,HM6,HM2}). In the present paper, there are used two different methods for construction of an almost contact B-metric manifold as a product of a real line and a 2-dimensional Norden manifold. The goal of this work is a characterisation of the obtained manifolds.

The paper is organized as follows. In
Sect.~\ref{sect-prel} we recall some preliminary facts about almost contact B-metric manifolds {and almost Norden manifolds}.
In Sect.~\ref{sect-C} and Sect.~\ref{sect-S} we study the considered manifolds derived from {a space-form} of real dimension 2 equipped with a complex structure and a Norden metric, which are constructed as {its cone} and {its }$S^1$-solvable extension, respectively.


\section{Preliminaries}\label{sect-prel}

\subsection{Almost contact B-metric manifolds}

Let $(M,\f,\xi,\eta,g)$ be an \emph{almost contact B-metric manifold}, where $M$
is a $(2n+1)$-dimensional differentiable manifold with an almost
contact structure $(\f,\xi,\eta)$ consisting of an endomorphism
$\f$ of the tangent bundle, a Reeb vector field $\xi$, its dual contact  1-form
$\eta$ as well as $M$ is equipped with a pseu\-do-Rie\-mannian
metric $g$ of signature $(n+1,n)$, called \emph{B-metric}, such that:
\begin{equation*}\label{strM}
\begin{array}{c}
\f\xi = 0,\qquad \f^2 = -\Id + \eta \otimes \xi,\qquad
\eta\circ\f=0,\qquad \eta(\xi)=1,\\[4pt]
g(\f x, \f y) = - g(x,y) + \eta(x)\eta(y),
\end{array}
\end{equation*}
where $\Id$ is the identity map (\cite{GaMiGr}). In the latter equality and further, $x$, $y$, $z$, $w$ will stand for arbitrary elements of $\X(M)$, the Lie algebra of tangent vector fields, or vectors in the tangent space $T_pM$ of $M$ at an arbitrary
point $p$ in $M$.
Let us recall that these manifolds are the odd dimensional extension of the almost Norden manifolds and the case with indefinite metrics corresponding
to almost contact metric manifolds.

The associated metric $\widetilde{g}$ of $g$ on $M$ is defined by
$\widetilde{g}(x,y)=g(x,\f y)+\eta(x)\eta(y)$ and {it is also a B-metric}. The manifold
$(M,\f,\xi,\eta,\widetilde{g})$ is also an almost contact B-metric
manifold.

The Ganchev-Mihova-Gribachev classification of almost contact B-metric manifolds, consisting of eleven basic classes $\F_1$, $\F_2$, $\dots$, $\F_{11}$, is given in
\cite{GaMiGr}. It is made with respect
to the (0,3)-tensor $F$ defined by
\begin{equation*}\label{F=nfi}
F(x,y,z)=g\bigl( \left( \n_x \f \right)y,z\bigr),
\end{equation*}
where $\n$ is the Levi-Civita connection of $g$ and the following {general properties are valid}:
\begin{equation*}\label{F-prop}
\begin{array}{l}
F(x,y,z)=F(x,z,y)=F(x,\f y,\f z)+\eta(y)F(x,\xi,z)
+\eta(z)F(x,y,\xi),\\[4pt]
F(x,\f y, \xi)=(\n_x\eta)(y)=g(\n_x\xi,y).
\end{array}
\end{equation*}

The intersection of the basic classes is the special class $\F_0$,
determined by the condition $F(x,y,z)=0$ and it is known as the
class of the \emph{cosymplectic B-metric manifolds}.

Let $\left\{e_i;\xi\right\}$ $(i=1,2,\dots,2n)$ be a basis of
$T_pM$ and let $\left(g_{ij}\right)$ be the {corresponding} matrix of $g$ and $\left(g^{ij}\right)$ be its inverse matrix. The 1-forms
$\theta$, $\theta^*$, $\omega$ {associated with $F$, called Lee forms, }are determined by:
\begin{equation*}\label{t}
\theta(z)=g^{ij}F(e_i,e_j,z),\quad \theta^*(z)=g^{ij}F(e_i,\f
e_j,z), \quad \omega(z)=F(\xi,\xi,z).
\end{equation*}

In the present work{,} we consider the case of the lowest dimension of the considered manifolds, \ie $\dim{M}=3$.

The basis $\left\{e_1,e_2,e_3\right\}$ in any tangent space at an arbitrary point of $M$ is called a \emph{$\f$-basis} if the following equalities are valid
\begin{equation}\label{fbasis}
\begin{array}{c}
\f e_1=e_2,\quad \f e_2=-e_1,\quad e_3=\xi,\\[4pt]
g(e_1,e_1)=-g(e_2,e_2)=g(e_3,e_3)=1,\quad
g(e_i,e_j)=0,\; i\neq j.
\end{array}
\end{equation}

According to \cite{HM1}, the components of $F$, $\ta$, $\ta^*$, $\om$, denoted by $F_{ijk}=F(e_i,e_j,e_k)$, $\ta_k=\ta(e_k)$, $\ta^*_k=\ta^*(e_k)$, $\om_k=\om(e_k)$,  with respect to the given $\f$-basis are:
\begin{equation}\label{t3}
\begin{array}{lll}
\ta_1=F_{221}-F_{331},\qquad &\ta_2=F_{222}-F_{332},\qquad &\ta_3=F_{223}-F_{322},\\[4pt]
\ta^*_1=F_{231}+F_{321},\qquad &\ta^*_2=F_{223}+F_{322},\qquad &\ta^*_3=F_{222}+F_{332},\\[4pt]
\om_1=0,\qquad &\om_2=F_{112},\qquad &\om_3=F_{113}.
\end{array}
\end{equation}

Let us denote by $F^s$ $(s=1,2,\dots,11)$ the components of $F$ in the
corresponding basic classes $\F_s$. {In \cite{HM1}, $F^s$ are obtained for $\dim M=3$ and it is established that }the class of 3-dimensional almost contact B-metric
manifolds is
\[
\F_1 \oplus \F_4 \oplus \F_5 \oplus \F_8 \oplus \F_9 \oplus
\F_{10} \oplus \F_{11}.
\]
{According to \cite{ManIv13}, the class of the normal almost contact B-metric manifolds is $\F_1\oplus\F_2\oplus\F_4\oplus\F_5\oplus\F_6$, since the Nijenhuis tensor of almost contact structure vanishes there. Therefore, we conclude that the class of the 3-dimensional normal almost contact B-metric manifolds is $\F_1\oplus\F_4\oplus\F_5$. The corresponding components of $F$ for $x=x^ie_i$, $y=y^je_j$, $z=z^ke_k$ are the following ones in compliance with \cite{HM1}:  }
\begin{equation}\label{Fi3}
\begin{array}{l}
F^{1}(x,y,z)=\left(x^2\ta_2-x^3\ta_3\right)\left(y^2z^2+y^3z^3\right),\\[4pt]
\qquad\ta_2=F_{222}=F_{233},\qquad \ta_3=-F_{322}=-F_{333}; \\[4pt]
F^{4}(x,y,z)=\frac{1}{2}\ta_1\Bigl\{x^2\left(y^1z^2+y^2z^1\right)
-x^3\left(y^1z^3+y^3z^1\right)\bigr\},\\[4pt]
\qquad \frac{1}{2}\ta_1=F_{212}=F_{221}=-F_{313}=-F_{331};\\[4pt]
F^{5}(x,y,z)=\frac{1}{2}\ta^*_1\bigl\{x^2\left(y^1z^3+y^3z^1\right)
+x^3\left(y^1z^2+y^2z^1\right)\bigr\},\\[4pt]
\qquad \frac{1}{2}\ta^*_1=F_{213}=F_{231}=F_{312}=F_{321}.
\end{array}
\end{equation}

In \cite{ManGri1}, there are given fundamental facts of the conformal geometry of almost contact B-metric manifolds. We recall some of them which are in relation with the obtained results in this work.

The contactly conformal transformations are determined by
\[
\bar{g}=e^{2u}\cos{2v}\ g+e^{2u}\sin{2v}\ \widetilde{g}+(1-e^{2u}\cos{2v}-e^{2u}\sin{2v})\eta\otimes\eta
\]
and they form a group $C$. The minimal contactly conformal equivalent class with respect to $C$ is contained in the most cramped class with non-zero Lee forms $\ta$ and $\ta^*$, the class $\F_1\oplus\F_4\oplus\F_5$, whose is also the minimal contactly conformal{ly} invariant class with respect to $C$.

The class{es} $\F_1\oplus\F_4$ and $\F_1\oplus\F_5$ are locally contactly conformally {invariant} with respect to the transformations of the groups $C_{1,4}$ and $C_{1,5}$, respectively, which are contactly conformal transformations satisfying the conditions $\D u(\xi)=0$ and $\D v(\xi)=0$, respectively. Similarly, the classes $\F_1$, $\F_4$ and $\F_5$ are locally contactly conformally {invariant} with respect to the transformations of the groups $C_1$, $C_4$ and $C_5$, respectively. The latter groups are determined by the conditions:
\begin{equation*}\label{contconf1,5}
\begin{array}{l}
C_1: \; du(\xi)=dv(\xi)=0,\\[4pt]
C_4: \; du\circ\f=dv\circ\f^2, \quad du(\xi)=0,\\[4pt]
C_5: \; du\circ\f=dv\circ\f^2, \quad dv(\xi)=0.
\end{array}
\end{equation*}
The subsets $\F_i^0$ ($i=1,4,5$) are subclasses of  $\F_i$ with closed Lee forms $\ta$ and $\ta^*$.
An $\F_0$-manifold is contactly conformally equivalent to a manifold from the classes $\F_1^0\oplus\F^0_4$, $\F_1^0\oplus\F^0_5$, $\F^0_1$, $\F^0_4$, $\F^0_5$ with respect to transformations belonging to the groups $C^0_{1,4}$, $C^0_{1,5}$, $C^0_1$, $C^0_4$, $C^0_5$, respectively.
The latter subgroups are determined as follows:
\begin{equation*}\label{contconf15}
\begin{array}{ll}
C_{1,4} \supset C^0_{1,4}: \; d(du\circ\f)=0, \quad &
C_{1,5} \supset C^0_{1,5}: \; d(dv\circ\f)=0, \\[4pt]
C_1 \supset C_1^0: \; d(du\circ\f)=d(dv\circ\f)=0,\quad \\[4pt]
C_4 \supset C_4^0: \; d(du\circ\f)=0,\quad &
C_5 \supset C_5^0: \; d(dv\circ\f)=0.
\end{array}
\end{equation*}

In \cite{Man31}, it is defined the square norm of $\n \f$ as follows
\begin{equation}\label{snf}
    \norm{\n \f}=g^{ij}g^{ks}
    g\bigl(\left(\n_{e_i} \f\right)e_k,\left(\n_{e_j}
    \f\right)e_s\bigr).
\end{equation}
If an almost contact B-metric manifold has a zero square
norm of $\n\f$ it is called an
\emph{isotropic-cosymplectic B-metric manifold} (\cite{Man31}).
Obviously, the equality $\norm{\n \f}=0$ is valid if $\M$
is an $\F_0$-manifold, but the inverse implication is not always
true.

%
%
%

Let $R=\left[\n,\n\right]-\n_{[\ ,\ ]}$ be the
curvature (1,3)-tensor of $\n$ and the corresponding curvature
$(0,4)$-tensor be denoted by the same letter: $R(x,y,z,w)$
$=g(R(x,y)z,w)$. The following properties are valid in general:
\begin{equation}\label{R}
\begin{array}{l}
R(x,y,z,w)=-R(y,x,z,w)=-R(x,y,w,z), \\[4pt]
R(x,y,z,w)+R(y,z,x,w)+R(z,x,y,w)=0.
\end{array}
\end{equation}

Let the essential curvature-like tensors $\pi_1$ and $\pi_2$ (\ie tensors generated by  $g$ and $\f$ such that they possess the properties \eqref{R}) of types (1,3) and (0,4) are defined by
\begin{equation}\label{pi12}
\begin{array}{l}
\pi_1(x,y)z=g(y,z)x-g(x,z)y,\\[4pt]
\pi_2(x,y)z=g(y,\f z)\f x-g(x,\f z)\f y,\\[4pt]
\pi_1(x,y,z,w)=g(y,z)g(x,w)-g(x,z)g(y,w),\\[4pt]
\pi_2(x,y,z,w)=g(y,\f z)g(x,\f w)-g(x,\f z)g(y,\f w).
\end{array}
\end{equation}

Let the Ricci
tensor $\rho$ and the scalar curvature $\tau$ for $R$ and $g$ as well as
their associated quantities be defined as follows
\begin{equation}\label{rhotau}
\begin{array}{c}
    \rho(y,z)=g^{ij}R(e_i,y,z,e_j),\qquad \rho^*(y,z)=g^{ij}R(e_i,y,z,\f e_j),\\[4pt]
    \tau=g^{ij}\rho(e_i,e_j),\quad
    \tau^*=g^{ij}\rho^*(e_i,e_j),\quad
    \tau^{**}=g^{ij}\rho^*(e_i,\f e_j).
\end{array}
\end{equation}

The sectional curvature of each non-degenerate 2-plane $\al$ in
$T_pM$ with respect to $g$ and $R$ has the following form
\begin{equation}\label{sect}
k(\al;p)=\frac{R(x,y,y,x)}{g(x,x)g(y,y)},
\end{equation}
where $\{x,y\}$ is an orthogonal basis of $\al$.

A 2-plane $\al$ is called a \emph{$\f$-holomorphic section}
(respectively, a \emph{$\xi$-section}) if $\al= \f\al$
(respectively, $\xi \in \al$).


\subsection{Almost complex manifold with Norden metric}

Let us remark that the $2n$-dimensional contact distribution
$H=\ker(\eta)$ of $(M,\f,\xi,\eta,g)$ can be considered as an almost complex manifold {$N$} endowed with an almost complex structure
$J=\f|_H$ and a metric $h=g|_H$, where $\f|_H$ and $g|_H$ are  the
restrictions of $\f$ and $g$ on $H$, respectively.

Let $x'$, $y'$, $z'$, $w'$ denote arbitrary vector fields or vectors in the contact distribution $H$ {of $M$}.
Since $g$ is a B-metric {of $M$}, then $h$ is a \emph{Norden metric} on $H$, \ie it
is compatible with $J$ as follows
\begin{equation*}\label{norden}
h(Jx',Jy')=-h(x',y').
\end{equation*}
The associated Norden metric $\widetilde{h}$ of $h$ is determined by
\begin{equation*}\label{assoc-norden}
\widetilde{h}(x',y')=h(x',Jy').
\end{equation*}
Both metrics $h$ and $\widetilde{h}$ have signature $(n,n)$.

We recall that an $2n$-dimensional manifold $N$ with almost complex structure $J$ and Norden metric $h$ is an almost Norden manifold $(N,J,h)$. The manifold $(N,J,\widetilde{h})$ is also an almost Norden manifold.


In the present work we pay attention to the case of the lowest dimension.
Let $(N,J,h)$ be a 2-dimensional almost Norden manifold. It is known that such a manifold is a space-form, \ie the manifold has constant sectional curvature $k'$ and its curvature tensor has the form $R'=k'\,\pi'_1$, where $\pi'_1$ is the essential curvature-like tensor as $\pi_1$ in \eqref{pi12} but with respect to $h$.
Moreover, according to \cite{GaBo}, any 2-dimensional almost Norden manifold is integrable and it belongs to the basic class $\W_1$.
In other words, this is the class of the considered manifolds which are conformally equivalent to K\"ahler-Norden manifolds, where conformal transformations of the metric are given by $\overline{h}=e^{2u}(\cos{2v}\,h+\sin{2v}\,\widetilde{h})$ for differentiable functions $u$, $v$ on $N$. The K\"ahler-Norden manifolds are the most specialized case of the almost Norden manifolds, determined by $\n' J=0$, where $\n'$ is the Levi-Civita connection of $h$. Their class is denoted by $\W_0$ and it is a subclass of $\W_1$. The fundamental tensor $F'$ of $(N,J,h)$ is defined by $F'(x',y',z')=h\bigl( \left( \n'_{x'} J \right)y',z'\bigr)$ and for any $\W_1$-manifold it is determined as follows
\begin{equation}\label{FW1}
\begin{array}{l}
  F'(x',y',z')=\frac12\{h(x',y')\ta'(z')+h(x',Jy')\ta'(Jz')\\[4pt]
  \phantom{F'(x',y',z')=\frac12\{}+h(x',z')\ta'(y')+h(x',Jz')\ta'(Jy')\}.
\end{array}
\end{equation}

{A} more wide subclass of $\W_1$ comparing with $\W_0$ contains the so-called \emph{isotropic-K\"ahler-Norden manifolds} (\cite{MekMan20}). They have a vanishing square norm of $\n' J$, \ie $\norm{\n' J}=0$,
where this
square norm is defined by
\begin{equation}\label{snJ}
  \norm{\n' J}=h^{ij}h^{ks}
    h\bigl((\n'_{e_i} J)e_k,(\n'_{e_j}J)e_s\bigr)
\end{equation}
with respect to an arbitrary basis.

\section{The cone over a 2-dimensional complex space-form with Norden metric}\label{sect-C}

In this section, we consider the cone over $N$, $\mathcal{C}(N)=\R^+\times N$, where $\R^+$ is the set of positive reals. We equip it with a metric $g$ defined by
\begin{equation}\label{C-g}
 g\left(\left(x',a\ddt\right),\left(y',b\ddt\right)\right)
=t^2\,h(x',y')+ab,
\end{equation}
where $t$ is the coordinate on $\R^+$ and $a$, $b$ are
differentiable functions on $\mathcal{C}(N)$.

We introduce an almost contact structure on the cone  by
\begin{equation}\label{C-str}
\f |_H=J, \quad \xi=\ddt, \quad \eta=\D t, \quad \f\xi=0,\quad \eta\circ\f =0.
\end{equation}
Obviously, the manifold $(\mathcal{C}(N), \f, \xi, \eta, g)$ is an almost contact B-metric manifold.

Using the general Koszul formula
\begin{equation}\label{koszul}
\begin{array}{l}
2g(\n_{x}y,z)=xg(y,z)+yg(z,x)-zg(x,y)\\[4pt]
\phantom{2g(\n_xy,z)=}+g([x,y],z)-g([y,z],x)+g([z,x],y),
\end{array}
\end{equation}
\eqref{C-g} and \eqref{C-str}, we obtain the following equalities for
the Levi-Civita connection $\n$ of the B-metric
$g$ on $\mathcal{C}(N)$:
\begin{equation}\label{C-n-g}
\begin{array}{ll}
    g\left(\n_{x'} y',z'\right)=t^2\,h\left(\n'_{x'} y',z'\right),\quad &
    g\left(\n_\xi y',z'\right)=t\,h\left(y', z'\right), %
    \\[4pt]
    g\left(\n_{x'} y',\xi\right)=-t\,h\left(x',y'\right),\quad &
    g\left(\n_{x'} \xi,z'\right)=t\,h\left(x', z'\right).
\end{array}
\end{equation}
Further, using \eqref{C-n-g}, we find the following formulae for the covariant derivatives with respect to $\n$
\begin{equation}\label{C-n}
    \n_{x'} y'=\n'_{x'} y'-\frac{1}{t}g\left(x', y'\right)\xi,\quad
    \n_\xi y'=\frac{1}{t}y',\quad
    \n_{x'} \xi=\frac{1}{t}x'.
\end{equation}

Bearing in mind \eqref{R} and \eqref{C-n}, we obtain
\begin{equation*}\label{C-Rxyz}
\begin{array}{l}
    R(x',y')z'=\frac{1}{t^2}(k'-1)\pi_1(x',y')z',\\[4pt]
    R(x',y')\xi=R(x',\xi)z'=R(\xi,y')z'=R(x',\xi)\xi=R(\xi,y')\xi=0.
\end{array}
\end{equation*}
Therefore, by direct computations, we get the following
\begin{prop}
The following equalities for the curvature tensor $R$ of $(\mathcal{C}(N), \f, \xi,\allowbreak{} \eta, g)$ are valid:
\begin{equation}\label{C-Rxyzw}
\begin{array}{l}
    R(x',y',z',w')=\frac{1}{t^2}(k'-1)\pi_1(x',y',z',w'),\\[4pt]
    R(\xi,y',z',w')=R(x',\xi,z',w')=R(x',y',\xi,w')\\[4pt]
    \phantom{R(\xi,y',z',w')}
    =R(x',y',z',\xi)=R(\xi,y',z',\xi)=0.
\end{array}
\end{equation}
\end{prop}

According to \eqref{fbasis} and \eqref{C-g}, we obtain the components $h_{ij}=h(e_i,e_j)$ and $g_{ij}=g(e_i,e_j)$ and the non-zero of them are
\begin{equation}\label{C-hijgij}
h_{11}=-h_{22}=\frac{1}{t^2}, \quad g_{11}=-g_{22}=g_{33}=1.
\end{equation}

Using \eqref{koszul}, \eqref{C-n} and \eqref{C-hijgij}, we get
\begin{equation}\label{C-nnij}
\begin{array}{lll}
\n_{e_1}e_1=\n'_{e_1}e_1-\frac{1}{t}e_3, \quad &\n_{e_1}e_2=\n'_{e_1}e_2, \quad &\n_{e_1}e_3=\frac{1}{t}e_1,\\[4pt]
\n_{e_2}e_1=\n'_{e_2}e_1, \quad &\n_{e_2}e_2=\n'_{e_2}e_2+\frac{1}{t}e_3, \quad &\n_{e_2}e_3=\frac{1}{t}e_2,\\[4pt]
\n_{e_3}e_1=\frac{1}{t}e_1, \quad &\n_{e_3}e_2=\frac{1}{t}e_2,  \quad &\n_{e_3}e_3=0.
\end{array}
\end{equation}

By virtue of \eqref{snf}, \eqref{C-str}, \eqref{C-hijgij} and \eqref{C-nnij}, we obtain the value of the square norm of $\n \f$ as follows
\begin{equation*}\label{C-nnf}
\begin{array}{l}
\norm{\n \f}=2\{(\ta_1)^2-(\ta_2)^2\}-\frac{4}{t^2},\\[4pt]
\norm{\n' J}=2t^2\{(\ta_1)^2-(\ta_2)^2\}.
\end{array}
\end{equation*}
The latter equalities imply
\begin{equation*}\label{C-nnfnnJ}
\norm{\n \f}=\frac{1}{t^2}\{\norm{\n' J}-4\}
\end{equation*}
and
the truthfulness of the following
\begin{thm}\label{C-nf-nJ}
\begin{enumerate}
  \item The manifold $(\mathcal{C}(N), \f, \xi, \eta, g)$ is an isotropic-cosymplectic B-metric manifold if and only if the square norm of $\n' J$ on $(N, J, h)$ is
\begin{equation*}
\norm{\n' J}=4.
\end{equation*}
  \item The manifold $(N, J, h)$ is an isotropic-K\"ahler-Norden manifold if and only if the square norm of $\n \f$ on $(\mathcal{C}(N), \f, \xi, \eta, g)$ is
\begin{equation*}
\norm{\n \f}=-\frac{4}{t^2}.
\end{equation*}
\end{enumerate}

\end{thm}

Taking into account \eqref{C-str}, \eqref{C-hijgij} and \eqref{C-nnij}, we compute the components $F_{ijk}$ of $F$ and get
\begin{equation}\label{C-Fijk}
\begin{array}{l}
  F_{111}=F_{122}=\ta'_1, \qquad  F_{211}=F_{222}=-\ta'_2,\\[4pt]
  F_{123}=F_{132}=F_{213}=F_{231}=\frac{1}{t},
\end{array}
\end{equation}
where $\ta'_i=\ta'(e_i)$ for $i=1,2$, as well as the other components of $F$ are zero.

Using \eqref{t3}, we obtain the components of the Lee forms of $(\mathcal{C}(N), \f, \xi,\allowbreak{} \eta, g)$ and the non-zero of them are:
\[
\ta_1=\ta'_1, \quad \ta_2=\ta'_2, \quad \ta^*_1=-\ta'_2, \quad \ta^*_2=\ta'_1, \quad \ta^*_3=\frac{2}{t}.
\]

Bearing in mind \eqref{Fi3} and \eqref{C-Fijk}, we establish the equality
\begin{equation*}\label{C-F1+F5}
F(x,y,z)=(F^1+F^5)(x,y,z),
\end{equation*}
where $F^1$ and $F^5$ are the components of $F$ in the basic classes $\F_1$ and $\F_5$, respectively.
The non-zero components of $F^1$ and $F^5$ by means of \eqref{Fi3} and \eqref{FW1} are the following
\begin{equation*}\label{Fijk15}
\begin{array}{l}
F^1_{111}=F^1_{122}=\ta_1, \quad F^1_{211}=F^1_{222}=-\ta_2,\\[4pt]
F^5_{123}=F^5_{132}=F^5_{213}=F^5_{231}=\frac{1}{2}\ta^*_3.
\end{array}
\end{equation*}
Therefore, we establish the truthfulness of the following
\begin{thm}
The manifold $(\mathcal{C}(N), \f, \xi, \eta, g)$
\begin{enumerate}
  \item belongs to $\F_1\oplus\F_5$,
  \item belongs to $\F_5$ if and only if $(N,J,h)$ is a $\W_0$-manifold,
  \item could not belongs to $\F_1$.
\end{enumerate}
\end{thm}

{The class $\F_1\oplus\F_5$ is a subclass of the class $\F_1\oplus\F_4\oplus\F_5$ of the 3-dimensional normal almost contact B-metric manifolds.}

Using \eqref{C-Rxyzw}, \eqref{C-hijgij} and \eqref{C-nnij}, we compute the components $R_{ijk\ell}=R(e_i,e_j,e_k,e_\ell)$ of the curvature tensor $R$. The non-zero ones of them are determined by \eqref{R} and the following
\begin{equation}\label{C-Rijkl}
  R_{1212}=\frac{1}{t^2}(k'-1).
\end{equation}

\begin{thm}
The manifold $(\mathcal{C}(N), \f, \xi, \eta, g)$ is flat if and only if $k'=1$.
\end{thm}

Bearing in mind \eqref{sect}, \eqref{C-hijgij} and \eqref{C-Rijkl}, we compute the basic sectional curvatures $k_{ij}=k(e_i,e_j)$ as follows
\begin{equation}\label{C-kij}
k_{12}=\frac{1}{t^2}(k'-1),\qquad k_{13}=k_{23}=0.
\end{equation}

Taking into account \eqref{rhotau}, \eqref{C-hijgij} and \eqref{C-Rijkl}, we obtain the basic components $\rho_{jk}=\rho(e_j,e_k)$ and $\rho^*_{jk}=\rho^*(e_j,e_k)$ of the
Ricci tensor $\rho$ and its associated tensor $\rho^*$, respectively, as well as the values
of the scalar curvature $\tau$ and its associated quantities $\tau^*$, $\tau^{**}$. The non-zero ones of them are:
\begin{equation}\label{C-rhotau}
\begin{array}{lll}
\rho_{11}=-\rho_{22}=\rho^{*}_{12}=\rho^{*}_{21}=\frac12\tau=\frac12\tau^{**}=\frac{1}{t^2}(k'-1).
\end{array}
\end{equation}

By virtue of \eqref{C-kij} and \eqref{C-rhotau}, we conclude the following
\begin{prop}
For the manifold $(\mathcal{C}(N), \f, \xi, \eta, g)$
\begin{enumerate}
  \item the sectional curvatures of the $\xi$-sections vanish,
  \item $\tau^*=0$,
  \item $\tau^{**}=\tau$.
\end{enumerate}
\end{prop}
\begin{prop}
The following assertions for $(\mathcal{C}(N), \f, \xi, \eta, g)$ are valid:
\begin{enumerate}
  \item $k'<1$ if and only if $\tau<0$;
  \item $k'=1$ if and only if $\tau=0$;
  \item $k'>1$ if and only if $\tau>0$.
\end{enumerate}
\end{prop}

\section{The $S^1$-solvable extension of a 2-dimensional complex space-form with Norden metric}\label{sect-S}

In the present section, let us introduce a warped product 3-dimensional manifold $S^1(N)=\mathbb R^+\times_{t^2} N$
as follows.
Let $\D t$ be the coordinate 1-form on $\mathbb R^+$ and let us define an
almost contact B-metric structure on $S^1(N)$ as follows
\begin{equation}\label{S-str}
\f |_H=J, \quad \xi=\ddt, \quad \eta=\D t, \quad \eta\circ\f =0,\quad g=\D t^2+\cos{2t}\,h-\sin{2t}\,\widetilde{h}.
\end{equation}

In \cite{IMM}, it is proved that the warped product manifold
$S^1(N)$ equipped with the almost contact B-metric structure defined in \eqref{S-str} is
an almost contact B-metric manifold. The constructed manifold $(S^1(N), \f, \xi, \eta, g)$ in this manner is called
\emph{an $S^1$-solvable extension of $(N,J,h)$} in \cite{IMM}.

Using \eqref{koszul} and \eqref{S-str}, we compute the components of the covariant derivative $\n$ as follows:
\begin{equation}\label{S-nn}
\begin{array}{c}
    \n_{x'} y'=\n'_{x'} y'- \frac{1}{2}\sin{2t}\{g\left(x',y'\right)\ta'^{\sharp}+ g\left(x',Jy'\right)J\ta'^{\sharp}\}+g\left(x',Jy'\right)\xi,\quad\\[4pt]
    \n_\xi y'=-Jy',\quad
    \n_{x'} \xi=-Jx',\quad \n_{\xi} \xi=0.
\end{array}
\end{equation}
 In the latter equalities and further, we denote by $\ta'^{\sharp}$ the dual vector of $\ta'$ with respect to $h$. Analogously, $\ta^{\sharp}$ stands for the dual vector of $\ta$ with respect to $g$.

Using \eqref{S-nn}, we obtain
\begin{equation}\label{S-Rxyz}
\begin{array}{l}
    R(x',y')z'=k'\{\cos2t\,\pi_1(x',y')z' - \sin2t\,J\pi_2(x',y')z'\} - \pi_2(x',y')z'\\[4pt]
    \phantom{R(x',y')z'=}+2\sin2t\,(\pi_1+\pi_2)(x',y',z',\ta'^{\sharp})J\ta^{\sharp}\\[4pt]
    \phantom{R(x',y')z'=}-\frac{1}{2}\sin{2t}\,\{g\left(y',z'\right)\n_{x'}\ta'^{\sharp}-g\left(x',z'\right)\n_{y'}\ta'^{\sharp}\\[4pt]
    \phantom{R(x',y')z'=k'\{\cos2t\,\pi}+g\left(y',Jz'\right)\n_{x'}J\ta'^{\sharp}-g\left(x',Jz'\right)\n_{y'}J\ta'^{\sharp}\}\\[4pt]
    \phantom{R(x',y')z'=}-\frac{1}{2}\{\cos2t\,(\pi_1+\pi_2)(x',y',z',\ta'^{\sharp})\\[4pt]
    \phantom{R(x',y')z'=k'\{}+\sin2t\,(\pi_1+\pi_2)(x',y',z',J\ta'^{\sharp})\}\xi,\\[4pt]
    R(x',y')\xi=-\frac{1}{2}(\pi'_1+\pi'_2)(x',y')\ta'^{\sharp},\\[4pt]
    R(\xi,y')z'=-\frac{1}{2}(\pi'_1+\pi'_2)(\ta'^{\sharp},y')z'\\[4pt]
    \phantom{R(\xi,y')z'=}-\cos2t\,\{h\left(y',z'\right)(\ta'_1e_1-\ta'_2e_2) + h\left(y',Jz'\right)(\ta'_2e_1+\ta'_1e_2)\}\\ \phantom{R(\xi,y')z'=}+g\left(y',z'\right)\xi,\\[4pt]
    R(x',\xi)z'=-\frac{1}{2}(\pi'_1+\pi'_2)(x',\ta'^{\sharp})z'\\[4pt]
    \phantom{R(x',\xi)z'=}+\cos2t\,\{h\left(x',z'\right)(\ta'_1e_1-\ta'_2e_2) + h\left(x',Jz'\right)(\ta'_2e_1+\ta'_1e_2)\}\\ \phantom{R(x',\xi)z'=}+g\left(x',z'\right)\xi.
\end{array}
\end{equation}

The equalities \eqref{S-str} and \eqref{S-Rxyz} imply the following
\begin{prop}
The following equalities for the curvature tensor $R$ of $(S^1(N), \f, \allowbreak{}\xi, \allowbreak{}\eta, g)$ are valid:
\begin{equation}\label{S-Rxyzw}
\begin{array}{l}
   R(x',y',z',w')=k'\{\cos2t\,\pi_1(x',y',z',w')-\sin2t\,\pi_2(x',y',z',\f w')\}\\[4pt]
    \phantom{R(x',y',z',w')=}-\pi_2(x',y',z',w')\\[4pt]
    \phantom{R(x',y',z',w')=}+4\sin2t\,\{\cos2t\,(\pi_1+\pi_2)(x',y',z',\ta^{\sharp})\\[4pt]
    \phantom{R(x',y',z',w')=k'\{\cos2t\,\pi}- \sin2t\,(\pi_1+\pi_2)(x',y',z',\f \ta^{\sharp})\}\ta(\f w')\\[4pt]
    \phantom{R(x',y',z',w')=}-\sin2t\,\{g\left(y',z'\right)[\cos2t\,(\n_{x'}\ta)w' + \sin2t\,(\n_{x'}\ta^*)w']\\[4pt]
    \phantom{R(x',y',z',w')=k'\{\cos2t}-g\left(x',z'\right)[\cos2t\,(\n_{y'}\ta)w' + \sin2t\,(\n_{y'}\ta^*)w']\\[4pt]
    \phantom{R(x',y',z',w')=k'\{\cos2t}+g\left(y',Jz'\right)[\sin2t\,(\n_{x'}\ta)w' - \cos2t\,(\n_{x'}\ta^*)w']\\[4pt]
    \phantom{R(x',y',z',w')=k'\{\cos2t}-g\left(x',Jz'\right)[\sin2t\,(\n_{y'}\ta)w' - \cos2t\,(\n_{y'}\ta^*)w']\},\\[4pt]
    R(x',y',z',\xi)=-(\pi_1+\pi_2)(x',y',z',\ta^{\sharp}),\\[4pt]
    R(x',y',\xi,w')=-(\pi_1+\pi_2)(x',y',\ta^{\sharp},w'),\\[4pt]
    R(x',\xi,z',w')=-(\pi_1+\pi_2)(x',\ta^{\sharp},z',w'),\\[4pt]
    R(\xi,y',z',w')=-(\pi_1+\pi_2)(\ta^{\sharp},y',z',w'),\\[4pt]
    R(\xi,y',z',\xi)=g\left(y',z'\right).
\end{array}
\end{equation}
\end{prop}


Using \eqref{fbasis} and \eqref{S-str}, we obtain the components $g_{ij}$ and $h_{ij}$ as follows
\begin{equation}\label{S-hijgij}
\begin{array}{ll}
g_{11}=-g_{22}=g_{33}=1, \quad &g_{12}=g_{21}=0,\\[4pt]
h_{11}=-h_{22}=\cos2t, \quad &h_{12}=h_{21}=-\sin2t.
\end{array}
\end{equation}

Taking into account \eqref{koszul}, \eqref{S-nn} and \eqref{S-hijgij}, we get
\begin{equation}\label{S-nnij}
\begin{array}{ll}
\n_{e_1}e_1=\n'_{e_1}e_1-\frac{1}{2}\sin2t\,\ta'^{\sharp}, \quad & \n_{e_1}e_2=\n'_{e_1}e_2+\frac{1}{2}\sin2t\,J\ta'^{\sharp}-\xi,\\[4pt]
\n_{e_2}e_1=\n'_{e_2}e_1+\frac{1}{2}\sin2t\,J\ta'^{\sharp}-\xi, \quad & \n_{e_2}e_2=\n'_{e_2}e_2+\frac{1}{2}\sin2t\,\ta'^{\sharp},\\[4pt]
\n_{e_1}e_3=\n_{e_3}e_1=-e_2, \quad & \n_{e_2}e_3=\n_{e_3}e_2=e_1,  \qquad \n_{e_3}e_3=0.
\end{array}
\end{equation}

Bearing in mind \eqref{snf}, \eqref{snJ}, \eqref{S-str}, \eqref{S-hijgij} and \eqref{S-nnij}, we obtain the values of the square norms of $\n \f$ and $\n' J$ as follows
\begin{equation*}\label{S-nnf}
\begin{array}{l}
\norm{\n \f}=2\{(\ta'_1)^2-(\ta'_2)^2\}+4, \\[4pt]
\norm{\n' J}=(1+\cos4t)\{(\ta'_1)^2-(\ta'_2)^2\}-2\sin4t\,\ta'_1\ta'_2.
\end{array}
\end{equation*}

The latter equalities imply the truthfulness of the following
\begin{thm}\label{nf-nJ}
 The manifold $(S^1(N), \f, \xi, \eta, g)$ is an isotropic-cosymplectic B-met\-ric manifold if and only if it is valid
\begin{equation*}
\norm{\n' J}=-2(1+\cos4t+\sin4t\,\ta'_1\ta'_2).
\end{equation*}
\end{thm}

By virtue of \eqref{S-str}, \eqref{S-hijgij} and \eqref{S-nnij}, we compute the components $F_{ijk}$ of $F$. The non-zero ones of them are
\begin{equation}\label{S-Fijk}
\begin{array}{l}
  F_{111}=F_{122}=\cos2t\,\ta'_1-\sin2t\,\ta'_2, \qquad  F_{211}=F_{222}=-\sin2t\,\ta'_1-\cos2t\,\ta'_2,\\[4pt]
  F_{131}=F_{113}=-F_{232}=-F_{223}=-1.
\end{array}
\end{equation}

Using \eqref{t3}, we obtain the components of the Lee forms of $(S^1(N), \f, \xi,\allowbreak{} \eta, g)$ and the non-zero of them are:
\[
\ta_1=\ta^*_2=\cos 2t\,\ta'_1-\sin 2t\,\ta'_2, \quad \ta_2=-\ta^*_1=\sin 2t\,\ta'_1+\cos 2t\,\ta'_2, \quad \ta_3=-2.
\]

Bearing in mind \eqref{Fi3} and \eqref{S-Fijk}, we establish the equality
\begin{equation*}\label{S-F1+F4}
F(x',y',z')=(F^1+F^4)(x',y',z'),
\end{equation*}
where $F^1$ and $F^4$ are the components of $F$ in the basic classes $\F_1$ and $\F_4$, respectively.
The non-zero components of $F^1$ and $F^4$ are
\begin{equation*}\label{Fijk14}
\begin{array}{l}
F^1_{111}=F^1_{122}=\ta_1, \quad F^1_{211}=F^1_{222}=-\ta_2,\\[4pt]
F^4_{131}=F^4_{113}=-F^4_{232}=-F^4_{223}=\frac{1}{2}\ta_3.
\end{array}
\end{equation*}
Therefore, we establish the validity of the following
\begin{thm}\label{thm14}
The manifold $(S^1(N), \f, \xi, \eta, g)$
\begin{enumerate}
  \item belongs to $\F_1\oplus\F_4$,
  \item belongs to $\F_4$ if and only if $(N,J,h)$ is a K\"ahler-Norden manifold,
  \item could not belongs to $\F_1$.
\end{enumerate}
\end{thm}

{The class $\F_1\oplus\F_4$ is a subclass of the class $\F_1\oplus\F_4\oplus\F_5$ of the 3-dimensional normal almost contact B-metric manifolds.}

Bearing in mind \eqref{S-Rxyzw}, \eqref{S-hijgij} and \eqref{S-nnij}, we calculate the components $R_{ijk\ell}$ of $R$. The non-zero of them are determined by \eqref{R} and the following
\begin{equation}\label{S-Rijkl}
\begin{array}{l}
  R_{1212}=k'\cos2t-1+\sin2t\,\cos2t\,\{(\n_{e_2}\ta)e_2-(\n_{e_1}\ta^*)e_2+8\ta_1\ta_2\}\\[4pt]
  \phantom{R_{1221}=k'\cos2t-1}+\sin^22t\,\{(\n_{e_2}\ta^*)e_2+(\n_{e_1}\ta)e_2+8(\ta_1)^2\},\\[4pt]
  R_{1213}=2\ta_2, \quad R_{1223}=2\ta_1, \quad R_{3113}=-R_{3223}=1.
\end{array}
\end{equation}

Bearing in mind \eqref{sect}, \eqref{S-hijgij} and \eqref{S-Rijkl}, we compute the basic sectional curvatures $k_{ij}$ as follows
\begin{equation}\label{S-kij}
k_{12}=R_{1212},\qquad k_{13}=k_{23}=1.
\end{equation}

Using \eqref{rhotau}, \eqref{S-hijgij} and \eqref{S-Rijkl}, we get the basic components $\rho_{jk}$ and $\rho^*_{jk}$ as well as the values $\tau$, $\tau^*$ and $\tau^{**}$:
\begin{equation}\label{S-rhotau}
\begin{array}{lll}
\rho_{11}=-\rho_{22}=R_{1212}+1, \quad & \rho_{33}=2, \quad &
\rho^{*}_{11}=\rho^{*}_{22}=\rho^{*}_{33}=0, \\[4pt]
\rho_{12}=\rho_{21}=0, \quad & \rho_{13}=\rho_{31}=-2\ta_1, \quad & \rho_{23}=\rho_{32}=-2\ta_2\\[4pt]
\rho^{*}_{12}=\rho^{*}_{21}=R_{1212}, \quad & \rho^{*}_{13}=\rho^{*}_{31}=2\ta_2, \quad & \rho^{*}_{23}=\rho^{*}_{32}=-2\ta_1,\\[4pt]
\tau=2R_{1212}+4, \quad & \tau^{*}=0, \quad & \tau^{**}=2R_{1212}.
\end{array}
\end{equation}

By virtue of the latter equalities for $\rho$, $\rho^*$ and \eqref{S-hijgij} for $g_{ij}$, we get
\begin{thm}
For the manifold $(S^1(N), \f, \xi, \eta, g)$ the following assertions are equivalent:
\begin{enumerate}
  \item $(N,J,h)$ is a K\"ahler-Norden manifold;
  \item $\rho=k'\cos2t\,g+(2-k'\cos2t)\eta\otimes\eta$;
  \item $\rho^*=(1-k'\cos2t)(\widetilde{g}-\eta\otimes\eta$).
\end{enumerate}
\end{thm}
Let us remark that the assertion (1) in the latter theorem is equivalent to the assertion that $(S^1(N), \f, \xi, \eta, g)$ is an $\F_4$-manifold, according to (2) in \thmref{thm14}.

{The assertion (2) in the latter theorem shows that $(S^1(N), \f, \xi, \eta, g)$ in this case is an $\eta$-Einstein manifold, according to \cite{HM2}; whereas the assertion (3) in the latter theorem presents $\rho^*$ as proportional to $g^*$, which is defined by $g^*(x,y)=g(x,\f y)$.}

By virtue of \eqref{S-kij} and \eqref{S-rhotau}, we conclude the following
\begin{prop}
For the manifold $(S^1(N), \f, \xi, \eta, g)$
\begin{enumerate}
  \item the sectional curvatures of the $\xi$-sections are constant,
  \item $\tau^*=0$,
  \item $\tau^{**}=\tau-4$.
\end{enumerate}
\end{prop}

Taking into account the first equality of \eqref{S-Rijkl} and the equalities in the last line of \eqref{S-rhotau}, we get the following corollaries.

\begin{cor}
If $(N,J,h)$ is a K\"ahler-Norden manifold, then
\[
\tau=2(k'\cos2t+1), \qquad \tau^{**}=2(k'\cos2t-1).
\]
\end{cor}

\begin{cor}
If $(N,J,h)$ is a K\"ahler-Norden manifold, then
\begin{enumerate}
  \item $k'<0$ if and only if $2(k'+1)<\tau<2$, $2(k'-1)<\tau^{**}<-2$;
  \item $k'=0$ if and only if $\tau=2$, $\tau^{**}=-2$;
  \item $k'>0$ if and only if $2<\tau<2(k'+1)$, $-2<\tau^{**}<2(k'-1)$.
\end{enumerate}
\end{cor}

\end{document}